\documentclass{article}
\usepackage{amsmath,graphicx}

\usepackage{amsthm}
\usepackage{amssymb}
\usepackage{xypic}
\usepackage{fancyhdr}
\usepackage{multirow}
\usepackage{hyperref}
\theoremstyle{plain}

\theoremstyle{definition}

\fancypagestyle{plain}{
\fancyhf{}
\lhead{}
\rhead{}
\chead{}
\lfoot{\tiny{This material is based upon work supported by the Defense Advanced Research Projects Agency (DARPA) under Contract No. HR001124C0319. Any opinions, findings and conclusions or recommendations expressed in this material are those of the author(s) and do not necessarily reflect the views of the Defense Advanced Research Projects Agency (DARPA). Distribution Statement ``A'' (Approved for Public Release, Distribution Unlimited).}}
}
\title{The structure of the token space for large language models}
\author{Michael Robinson$^1$ \and Sourya Dey$^2$ \and Shauna Sweet$^3$}
\date{%
  $^1$Mathematics and Statistics, American University, Washington, DC, USA, michaelr@american.edu\\%
  $^2$Galois, Inc., Arlington, VA, USA, sourya@galois.com\\%
  $^3$sjaynesweet@gmail.com\\%
}

\begin{document}

\maketitle

\begin{abstract}
Large language models encode the correlational structure present in natural language by fitting segments of utterances (\emph{tokens}) into a high dimensional ambient \emph{latent space}
upon which the models then operate.
We assert that in order to develop a foundational, first-principles understanding of the behavior and limitations of large language models,
it is crucial to understand the topological and geometric structure of this \emph{token subspace}.
In this article, we present estimators for the dimension and Ricci scalar curvature of the token subspace, and apply it to  three open source large language models of moderate size: GPT2, LLEMMA7B, and MISTRAL7B.
In all three models, using these measurements, we find that the token subspace is not a manifold, but is instead a stratified manifold,
where on each of the individual strata, the Ricci curvature is significantly negative.
We additionally find that the dimension and curvature correlate with generative fluency of the models, which suggest that these findings have implications for model behavior.
\end{abstract}

\pagestyle{plain}

\tableofcontents

\section{Introduction}

Large language models (LLMs) encode the correlational structure present in natural language by fitting segments of utterances (\emph{tokens}) into a high dimensional ambient \emph{latent space}
upon which the models then operate.
Not every point in this latent space is linguistically meaningful; only a subspace of it corresponds to tokens that are extant in the learned vocabulary of the target language.
We assert that in order to develop a foundational, first-principles understanding of the behavior and limitations of large language models,
it is crucial to understand the topological and geometric structure of this \emph{token subspace}.

Since the tokens are determined prior to model training,
one of the tasks of training from scratch is to determine numerical coordinates for each token.
Operationally, the token subspace is therefore defined via the image of a learned \emph{embedding function} (or simply an \emph{embedding}), which defines the values of the numerical coordinates at each token.
In this paper, we take this embedding function as our starting point.

The mathematical definition of an embedding requires that both the ambient latent space and the space of tokens each have a well-defined \emph{topology}.
Since the latent space is usually ascribed the usual Euclidean metric, this induces a topology on both the latent space and the token subspace.
It is therefore most reasonable to assert that the space of tokens, abstractly,
also has a topology induced by the embedding function, which satisfies the mathematical definition of an embedding.
From these assertions, we can then show that it is possible to reliably estimate two properties of the token subspace: \emph{dimension} (a topological property) and \emph{Ricci scalar curvature} (a geometric property). 
We aim to show that with these estimates of dimension and Ricci scalar curvature, it is possible to reliably anticipate the inferential quality of model behavior within and across regions of the token subspace.

In this article, we present the dimension and Ricci scalar curvature for three open source large language models of moderate size: GPT2 \cite{radford2019language}, LLEMMA7B \cite{azerbayev2023llemma}, and MISTRAL7B \cite{jiang2023mistral}.
In all three models, using these measurements, we find that \emph{the token subspace is not a manifold, but is instead a stratified manifold},
where on each of the individual strata, \emph{the Ricci curvature is significantly negative.}
These findings have implications for model behavior, its boundedness and predictability at inference, and its stability under retraining.

\subsection{Contributions}

It is a natural question whether the token subspace has the structure of a \emph{manifold},
a topological space for which every point has a neighborhood that is topologically equivalent to Euclidean space \cite{Lee_2003}.
\emph{We demonstrate instead that the space of tokens is not a manifold for several large language models.}
The implications are striking: some tokens have dramatically more non-token neighbors in the latent space than do other tokens,
and this changes abruptly as one ``moves around'' within the space of tokens.
Since the transformer blocks that are used to predict subsequent tokens employ continuous functions \cite{vaswani2023attentionneed},
these abrupt changes in dimension imply that the behavior of the large language model as a whole is unstable to small perturbations.
Because the space of tokens is a stratified manifold, model behavior in response to queries is necessarily locally, not globally, determined.

We also find \emph{that the dimension of the token subspace, estimated at each token, is much lower than the dimension of the ambient latent space}.
Moreover, we find that the dimension of the token subspace is highly variable in ways that are correlated with generative fluency.
In this article, we examine the numeric versus non-numeric tokens; across models these have substantially different local dimension.
In GPT2, the numeric tokens are confined to a low dimensional stratum (an embedded submanifold of the latent space),
whereas in LLEMMA7B and MISTRAL7B, the numeric tokens do not form a connected subspace, and moreover nearly half of them are isolated from all other tokens.

In addition to dimension, which is a topological property, the token subspace has estimable geometric properties, which are measurable and also have implications for model behavior.
Since the Euclidean metric defines distances and angles, it can likewise induce geometric information on the space of tokens.
The most easy property to assess is curvature, which like dimension, is related to the size of neighborhoods of a token.
\emph{This article shows that the Ricci scalar curvature is always significantly negative for several large language models, yet substantial variations in curvature exist across strata.}
Again, there is a correlation between curvature and generative fluency.
In GPT2, the stratum along which numeric tokens are concentrated is significantly flatter (curvature nearer to zero) than the other strata,
whereas in MISTRAL7B and LLEMMA7B, the numeric tokens are associated with roughly the same curvature as the other tokens.


\subsection{Related research}

Metrics for language have been considered before large language models became popular \cite{Stratford_2009}.
The problem of finding the token embedding itself is a special case of matrix factorization \cite{Omer_2015,Omer_2014},
in which the subspace of tokens is assumed to span a low-dimensional linear subspace, which is a manifold \cite{bradley2023structure}.
Aside from \cite{Kohn_2024,jakubowski_2020}, in which it is argued that subspace of neuron activations or tokens might be a manifold with singularities, we are unaware of any work that attempts to
interrogate the assumption that the relevant subspace of the latent space is a manifold, or attempts to discover the topological or geometric structure of the token subspace itself.
This is surprising since a hypothesis test for manifolds is known \cite{Fefferman_2016},
and the dimension of word embeddings has been estimated before (for instance \cite{Gromov_2024,tulchinskii2023intrinsicdimensionestimationrobust}, among others).
This work both attempts to detect and characterize these structures, as well as link them to model behavior. 

The idea that not every point in the ambient latent space corresponds to a token is consistent in what is observed in the image domain.
The paper \cite{tu2023probabilistic} notes that ``not every pattern of pixels is an image. It is common to say that images lie on a lower-dimensional manifold in the high dimensional space.''
Our work highlights what might be a critical distinction with how semantics are learned by transformers,
namely that tokens do not lie on a lower-dimensional manifold,
but instead on a stratified manifold.

While it seems reasonable to use topological methods, such as persistent homology to estimate dimension, as is done in \cite{Rohrscheidt_2022, jakubowski_2020, schweinhart2020fractaldimensionpersistenthomology}, these methods tend to scale poorly.
The reason for scaling difficulty is that generally one must build a representation of the topological space,
which involves constructing myriad simplices.
Instead, we rely on a regression-based method that does not require this construction, and so is practical enough to be deployed at scale.

\subsection{Implications}

If the token subspace is not a manifold, this has important implications because the behavior of the transformer blocks\footnote{We note that although this paper examines models that use transformers to generate textual or linguistic data, similar architecture components are used in models that generate images.},
which are piecewise smooth (hence continuous) transformations of the latent space \cite{villani_2024},
must therefore preserve the dimensions we observe. 
As a result, queries that cross stratification boundaries will yield responses that exhibit dramatic changes in behavior.
This instability will likely \emph{preclude strong guarantees about the model's generative performance} without intimate knowledge of how the token subspace is embedded within the ambient latent space.

\section{Methods}

In this section, we describe a novel Monte-Carlo-based method for estimating both local dimension and Ricci scalar curvature of the token subspace.
We also briefly summarize approaches to interpret estimates using both visual and parametric methods.

First, we introduce the relationship between volume and radius as a function of dimension and curvature, our specific quantities of interest.
In Euclidean space of dimension $d$, each point corresponds to a vector in $\mathbb{R}^d$, in which distances between two points $x, y \in \mathbb{R}^d$ is given by the familiar formula
\begin{equation*}
  dist(x,y) = \sqrt{\sum_{i=1}^d |x_i - y_i|^2}.
\end{equation*}
In this setting, the volume $v$ of a ball of radius $r$ is
\begin{equation}
  \label{eq:volume}
  v = \frac{\pi^{d/2}}{\Gamma(d/2+1)}r^d.
\end{equation}
It follows that one can estimate $d$ from comparing estimates of volume with radius.

If the token subspace is not Euclidean, but is still a manifold, the radial dependence of volume changes.
Moreover if the token subspace is a stratified manifold, in the limit of small $r$, the volume at any point in the interior of a stratum\footnote{In a stratified manifold, the points on stratification boundaries form a set of measure zero.} is asymptotic to (see \cite[Thm 3.1]{Gray_1974} if the space is manifold, \cite[Cor. 5.2]{Yomdin_2004} otherwise)
\begin{equation}
  \label{eq:volume_curvature}
  \begin{aligned}
    v =  K r^n \left(1 - \frac{1}{6(n+2)} Ric \; r^2 + O(r^4) \right),
  \end{aligned}
\end{equation}
where $n$ is the local dimension of the token subspace, $Ric$ is the Ricci scalar curvature of the token subspace, and $K$ is constant of proportionality.
We will call $K$ the \emph{volume scaling coefficient}, noting here that its value is difficult to estimate without prior knowledge of the space's total volume.

Taking the natural logarithm of both sides of Equation \eqref{eq:volume_curvature} yields the following asymptotic series for small $r$,
\begin{equation}
  \label{eq:volume_curvature_log}
  \log v = \log K + n \log r - \frac{Ric}{6(n+2)} r^2 + O(r^4),
\end{equation}
where $\log v$ depends linearly on the three desired parameters ($\log K$, $n$, and $Ric$),
we may solve for each of these parameters via a linear regression against pairs of radius-versus-volume values.
Recognizing the potential numerical instability of regression estimates, due to the dynamic range in the term involving the Ricci scalar curvature,
we solve for each of these parameters as follows:
first determining the scaling coefficient $K$ and the dimension $n$ by solving a linear regression,
and then from the residual determining $Ric$.

We correct for bias in $K$ using the procedure in \cite{Miller_1984},
\begin{equation*}
  \widehat{K}' = \widehat{K} e^{\sigma^2/2},
\end{equation*}
where $\sigma$ is the standard error of the intercept in the regression.

The residual for the regression problem for the first two terms on the right in Equation \eqref{eq:volume_curvature_log} will then be dominated by the next term in the series, namely the term involving the Ricci scalar curvature.
That is, if $\widehat{n}$ and $\widehat{K}'$ are the estimates obtained by regression,
we can therefore estimate the Ricci scalar curvature by
\begin{equation}
  \label{eq:ricci_estimator}
  \widehat{Ric} = \text{mean } \left\{\frac{6(\widehat{n} + 2)}{r^2} \left(\log \widehat{K}' + \widehat{n} \log r - \log v\right)\right\}.
\end{equation}

\subsection{Estimating volume}
\label{sec:alg}

Equation \eqref{eq:volume_curvature_log} gives a powerful tool for extracting topological and geometric features from a log-log plot of the volume versus radius of a disk centered at a particular token within the token subspace.
The token embedding produces a discrete subset of Euclidean space (a so called \emph{point cloud}),
which has the appearance of being a set of samples drawn from a larger subspace that ``interpolates them'' in the latent space.
If we assert that the tokens are approximately uniformly distributed on the token subspace\footnote{If the tokens are not in general position, the strict inequalities may degenerate to equalities for some values of $k$.  This poses no problem for our method.},
then a Monte Carlo estimation of volume is reasonable.
The volume of a ball of radius $r$ centered at a token $j$ is proportional to the number of points within a distance $r$ to $j$,
that is
\begin{equation*}
  v(r;j) \approx M \#\{ i : \|i-j\| \le r\},
\end{equation*}
where the vertical bars indicate a distance calculation, and $M$ is the volume contribution of each token to the Monte-Carlo estimate.
Except in the case of subspaces of known volume, $M$ is usually not known and must be treated as a nuisance constant.

Suppose that there are $p$ tokens in total and that the tokens are located in general position in the latent space.
Then for a given token $j$, we can obtain a sequence
\begin{equation*}
  r_{1,j} < r_{2,j} < \dotsb < r_{p,j}
\end{equation*}
of distances to the other tokens such that
\begin{equation*}
  v(r_{k,j}; j) \approx M k.
\end{equation*}

We can arrange the set of $r_{i,j}$ values in a matrix, in which we interpret $i$ as specifying rows and $j$ as specifying columns.
Notice that each column of the $r_{i,j}$ matrix is sorted in ascending order.
Since the sequence of volumes corresponding each token (column) is the same, namely the sequence of integers from $1$ to $p$,
we can rewrite the log-linear portion of Equation \eqref{eq:volume_curvature_log} as a matrix equation for token $j$,
\begin{equation}
  \label{eq:volume_log_matrix}
  \begin{pmatrix}
    0 \\ \log 2 \\ \vdots \\ \log p \\
  \end{pmatrix} + \log M
  \approx
  \begin{pmatrix}
    1 & \log r_{1,j} \\
    1 & \log r_{2,j} \\
    \vdots & \vdots \\
    1 & \log r_{p,j} \\
  \end{pmatrix}
  \begin{pmatrix}
    \log \widehat{K_j} \\
    \widehat{n_j}
  \end{pmatrix}
  + O(r^2)
\end{equation}
This is readily solved via least squares regression to obtain estimates for the dimension $\widehat{n_j}$ and scaling coefficient $\widehat{K_j}$ at each token $j$.

If $M$ is unknown, Equation \eqref{eq:volume_log_matrix} makes it clear that the estimates of $\widehat{K_j}$ are determined up to a multiplicative constant, and moreover $\widehat{n_j}$ is not impacted.
The estimate of $\widehat{Ric}$ in Equation \eqref{eq:ricci_estimator} is also not impacted because it only depends upon the difference $\log \widehat{K_j} - \log v$, in which the contributions of $M$ cancel.

To reduce the impact of sampling error from balls of small radius, and to reduce the impact of curvature and stratifications on balls of large radius,
it is advantageous to solve the regression problem using only a band from the middle rows in Equation \ref{eq:volume_log_matrix}.
We found during this work that the precise band of volumes to be used depends delicately upon the particular volume-versus-radius curves,
with the starting row adjusted to avoid contamination from points of different strata.

\subsection{Interpreting volume-versus-radius curves}
\label{sec:v_r}

The log-log plot of volume versus radius at each token provides considerable information.
Briefly,
\begin{enumerate}
\item The slope of the curve corresponds to local dimension,
\item A ``knee'' in the curve, where the slope changes from one non-negative value to another is a clear sign that the space is not a manifold,
\item A horizontal gap in the curve (a zero slope portion) may mean that the neighborhood of the token contains multiple connected components, and
\item The concavity of the curve determines the local curvature, where concave up corresponds to negative Ricci scalar curvature, and concave down corresponds to positive Ricci scalar curvature.
\end{enumerate}

Unpacking each of these in turn, it is immediately clear that the slope (derivative) of such a curve for small $r$ (near the bottom left of the plot), gives the dimension of the neighborhood of the token.
Since our estimates are obtained by Monte-Carlo estimation, it is wise to exclude a small neighborhood of the token from analysis, since this tends to have large sampling error.

If the slope changes, this means that the dimension of the neighorhood of the token has changed as the radius increases.
It is possible for the slope to increase or decrease, depending on whether the neighborhood expands to contain a greater or lesser dimension region.
It is worth noting that estimates of dimension far from the original center token may be biased, but the fact of an abrupt change is compelling evidence that the token subspace is not a manifold.

On the other hand, a gap in the token subspace will mean that there are no additional tokens within a certain band of radii.
This means that the volume estimate will not change over these radii.
The result is a curve with horizontal (zero) slope.

Finally, according to Equation \eqref{eq:volume_curvature_log}, nonzero Ricci curvature adds a quadratic term to the volume.
This makes the volume-versus-radius curve either concave up or down according to the sign of the curvature.

\section{Results}

\begin{table}
  \caption{Parameters of LLMs under test}
  \label{tab:llm_parameters}
  \begin{tabular}{|l|c|c|c||c|}
    \hline
    Model & Latent dim. & \multicolumn{3}{|c|}{Token counts}\\
    &             & Non-numeric & Numeric & Total \\
    \hline
    \hline
    GPT2 & 768 & 48563 & 1694 & 50257\\ 
    LLEMMA7B & 4096 & 31972 & 44 & 32016 \\ 
    MISTRAL7B & 4096 & 31985 & 31 & 32016 \\
    \hline
  \end{tabular}
\end{table}

We examined the token subspaces of three LLMs: GPT2, LLEMMA7B, and MISTRAL7B.
The basic features of these subspaces are shown in Table \ref{tab:llm_parameters}.
A token that contains any numeric character is marked as a ``numeric token'', otherwise it is marked as a ``non-numeric token''.
There are very few numeric tokens in LLEMMA7B and MISTRAL7B,
so the reader is cautioned that statistical inferences regarding numeric tokens are subject to sampling error.

\begin{table}
  \caption{Estimated parameters of token space}
  \label{tab:llm_token_summary}
  \begin{tabular}{|l|c||c|c|c||c|c|c|}
    \hline
    Model & Subset & \multicolumn{3}{|c||}{Dimension} & \multicolumn{3}{|c|}{Ricci scalar curvature}  \\
    &           & Q1 & Q2 & Q3 & Q1 & Q2 & Q3\\
    \hline
    \hline
    GPT2 & Non-numeric  & 384 & 498 & 566 & -107 & -63 &  -31.7  \\
    & Numeric & 10.7 & 15.8 & 22.1 & -4.68 & -2.39 & -1.40 \\
    \hline

    LLEMMA7B & Non-numeric & 9.44 & 10.5 & 11.2 &     -185 & -169 & -153 \\
    & Numeric & 4.92 & 6.84 & 9.07    & -268 & -170 & -154\\
    \hline

    MISTRAL7B & Non-numeric & 5.14 & 5.64 & 6.07 & -5339 & -5036 & -4753 \\
    & Numeric & 0.510 & 2.83 & 5.21 & -82938 & -5693 & -4782\\
    \hline
    \end{tabular}
\end{table}

A summary of our estimates of dimension and Ricci scalar curvatures are shown in Table \ref{tab:llm_token_summary}.
Because the total volumes of each token subspace are not known in advance, we cannot estimate the scaling coefficients.

The distribution of local dimensions of the three models are very statistically significantly different ($p <0.001$ for the Kolmogorov-Smirnov test).
In LLEMMA7B and MISTRAL7B the local dimensions of the token subspace are well below half of the latent space dimension,
whereas for GPT2 the typical dimension is higher.
It is worth noting that the estimates for LLEMMA7B and MISTRAL7B are broadly consistent with the results obtained from other word embeddings \cite{Gromov_2024,tulchinskii2023intrinsicdimensionestimationrobust}.
In all of these cases, the conclusion is the same: most of the latent space is not ``near'' a token in any meaningful way.
From a practical standpoint, this means that the embeddings are quite wasteful of memory, especially for LLEMMA7B and MISTRAL7B.

For GPT2, the Kolmogorov-Smirnov test indicates that the distribution of local dimensions for the numeric and non-numeric tokens differ (with $p<0.001$).
(Given the small number of numeric tokens in LLEMMA7B and MISTRAL7B, while the dimensions of numeric and non-numeric tokens are statistically significant to the same level, any such comparison is subject to substantial sampling error.)
Moreover, as is discussed in the subsections below, the numeric tokens of GPT2 are localized to a small number of connected components with relative uniform dimension.
Notice that the interquartile range of dimension for GPT2 numeric tokens is $11.4$, whereas for non-numeric tokens it is nearly $200$.
In contrast, the numeric tokens of LLEMMA7B and MISTRAL7B are frequently isolated points within the token subspace (see Figures \ref{fig:llemma7b_dim_hist} and \ref{fig:mistral7b_dim_hist}).

As the subsections below show, inspection of volume-versus-radius curves show that the local dimension varies significantly over connected components of the token subspaces in all three LLMs.
This gives strong evidence that the token subspaces are not manifolds, but are stratified manifolds.

\subsection{GPT2}

\begin{figure}[!htbp]
  \begin{center}
    \includegraphics[width=3in]{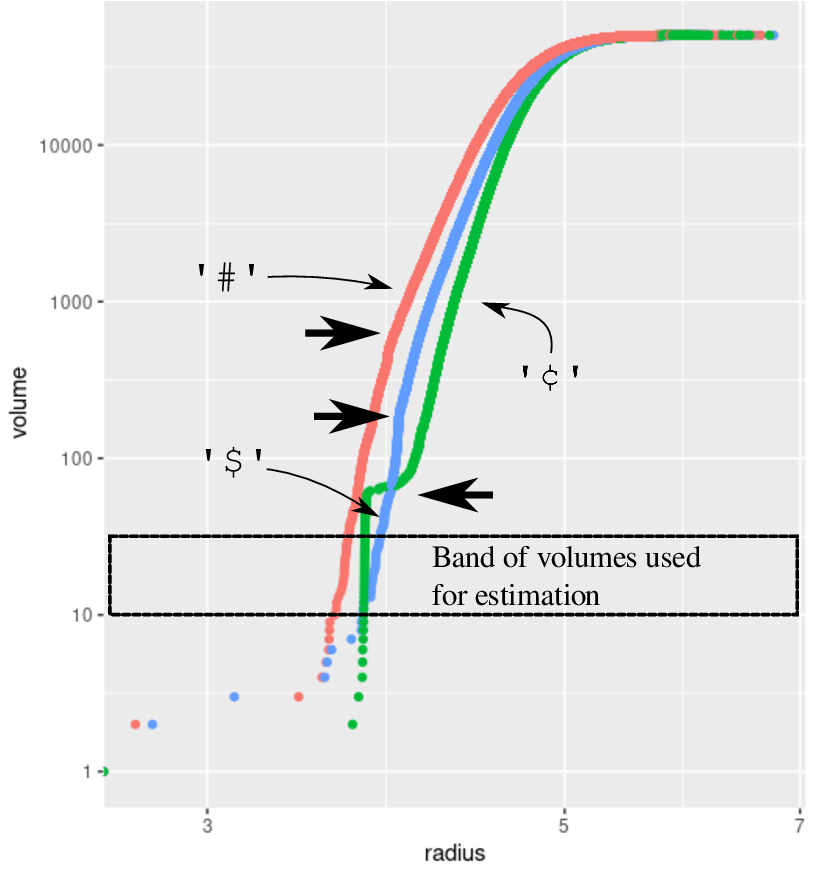}
    \caption{Volume-versus-radius for several tokens in GPT2: {\tt \#} (red), {\tt \$} (blue), {\tt \textcent} (green).  Visible stratification boundaries are marked with arrows.}
    \label{fig:gpt2_v_r}
  \end{center}
\end{figure}

Figure \ref{fig:gpt2_v_r} shows volume-versus-radius curves for several tokens in GPT2.
Because of Monte-Carlo estimation error, there is some variability for small volumes,
so we begin our estimation with the volume corresponding to $10$ tokens to avoid undue error.

Knees between mostly straight segments are prominently visible in the curves for all tokens shown,
which provides compelling evidence that the token subspace is a stratified manifold.
It is also interesting that the dimension of the tokens, and the structure of their stratifications, is different for different tokens.
For instance, the difference between the US dollar symbol {\tt \$} versus other currency symbols may be due to its use in the syntax of programming languages.

\begin{figure}[!htbp]
  \begin{center}
    \includegraphics[width=3in]{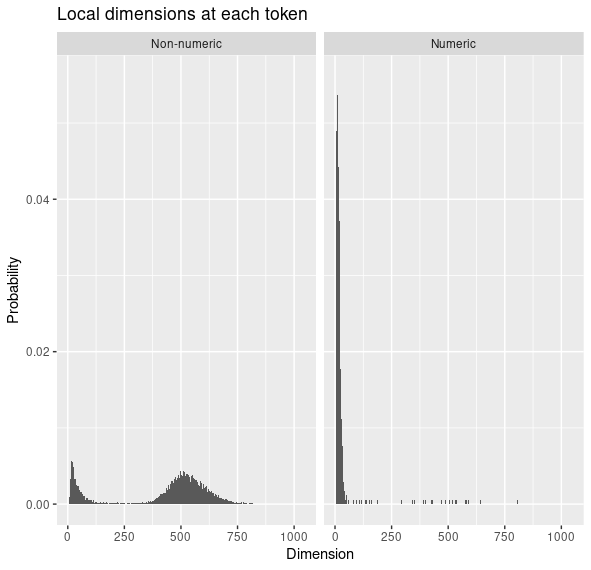}
    \caption{Histogram of estimated local dimensions for GPT2.}
    \label{fig:gpt2_dim_hist}
  \end{center}
\end{figure}

Figure \ref{fig:gpt2_dim_hist} shows the distribution of local dimensions estimated by our method for GPT2's token subspace.
The distribution of numeric tokens is unimodal and skewed rightward,
while the distribution of non-numeric tokens is clearly bimodal.
Interestingly, the cluster of non-numeric tokens with low dimension largely consists of tokens related to dates and current events, and mostly consists of whole words.
Tokens that are word fragments appear to become increasingly common as the dimension increases.

\begin{figure}[!htbp]
  \begin{center}
    \includegraphics[width=4in]{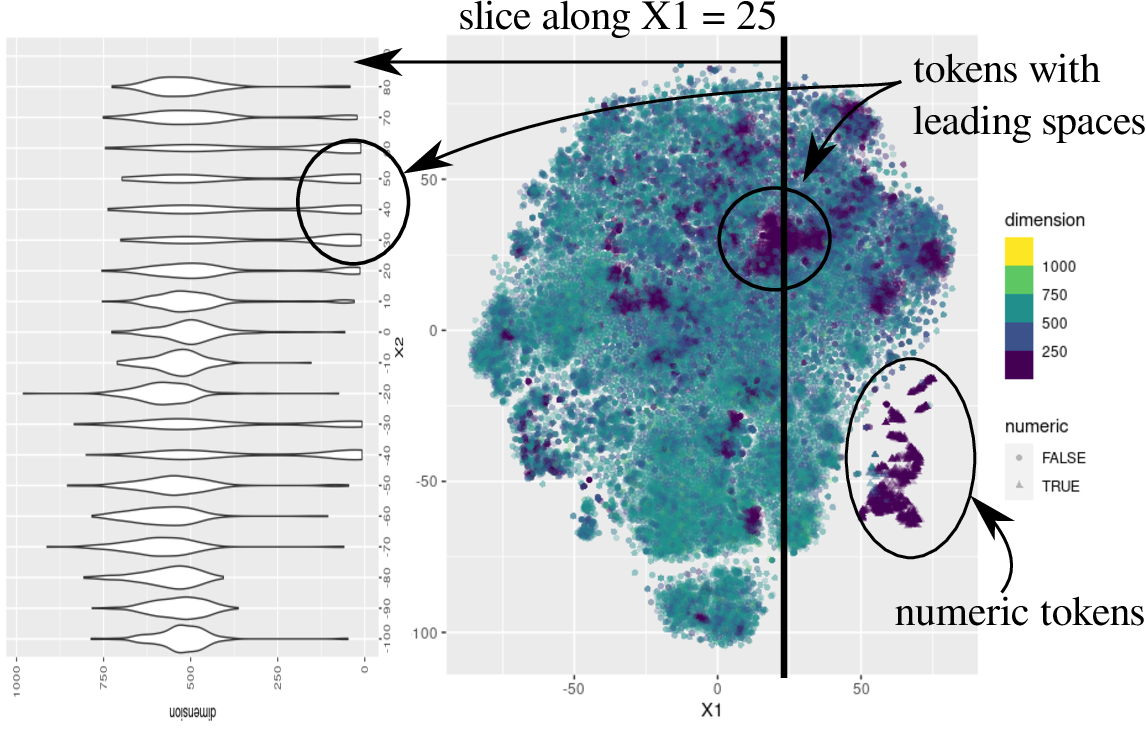}
    \caption{Estimated local dimension for GPT2, plotted using tSNE coordinates.  Dimension is indicated by the depth of color: darker points have lower local dimension.  Numeric tokens are all in the region marked.}
    \label{fig:gpt2_dim_overview}
  \end{center}
\end{figure}

Because the local dimension of the token subspace is well above $2$, it is impossible to render the token subspace isometrically in $2$ dimensions.
In order to obtain a visualization, some concessions must be made that will introduce distortions.
Figure \ref{fig:gpt2_dim_overview} shows a tSNE projection\cite{Maaten_2008} of the token subspace, colored by local dimension.
Each point in the figure is a distinct token.
Since tSNE is a continuous transformation, any tokens that are nearby each other in the token subspace will be nearby their images in Figure \ref{fig:gpt2_dim_overview}.
The converse does not hold; there can be tokens whose corresponding images in Figure \ref{fig:gpt2_dim_overview} are nearby each other, but whose actual distance apart in the token subspace is large.
However, the usual concerns with clustering and tSNE are that it distorts the geometry of and between connected components,
sometimes to the point of forming spurious connected components \cite{Chari_2023,Schubert_2017}.
In our case, we are more concerned with the opposite problem, that tSNE might combine connected components.
While this is possible with arbitrary continuous maps, it is substantially less likely with tSNE itself,
because separate connected components will tend to repel each other in the tSNE cost function.

Most of the numeric tokens fall in the region indicated, though in that region there are also non-numeric tokens present, and a few numeric tokens are scattered elsewhere in the token space.
That region mostly consists of a single connected component in Figure \ref{fig:gpt2_dim_overview}.
Continuity of the tSNE transformation does not guarantee that this region is a single connected component,
but because of the way tSNE is actually constructed, it is extremely unlikely that there are multiple connnected components.

Figure \ref{fig:gpt2_v_r} provides conclusive proof that the token subspace cannot be a manifold.
This is also confirmed by examining a vertical slice of Figure \ref{fig:gpt2_dim_overview}, and aggregating the nearby local dimensions along this slice.
Figure \ref{fig:gpt2_dim_overview} shows how the distribution of local dimensions along this slice vary, when grouped into bins of size $10$ along the $X1=25$ tSNE axis.
Notice that there is an abrupt, and statistically significant, change in local dimension that happens near $X2=50$.
This strongly suggests that the token subspace does not have a single local dimension within the corresponding connected component, and therefore is not a manifold.

\begin{figure}[!htbp]
  \begin{center}
    \includegraphics[width=3in]{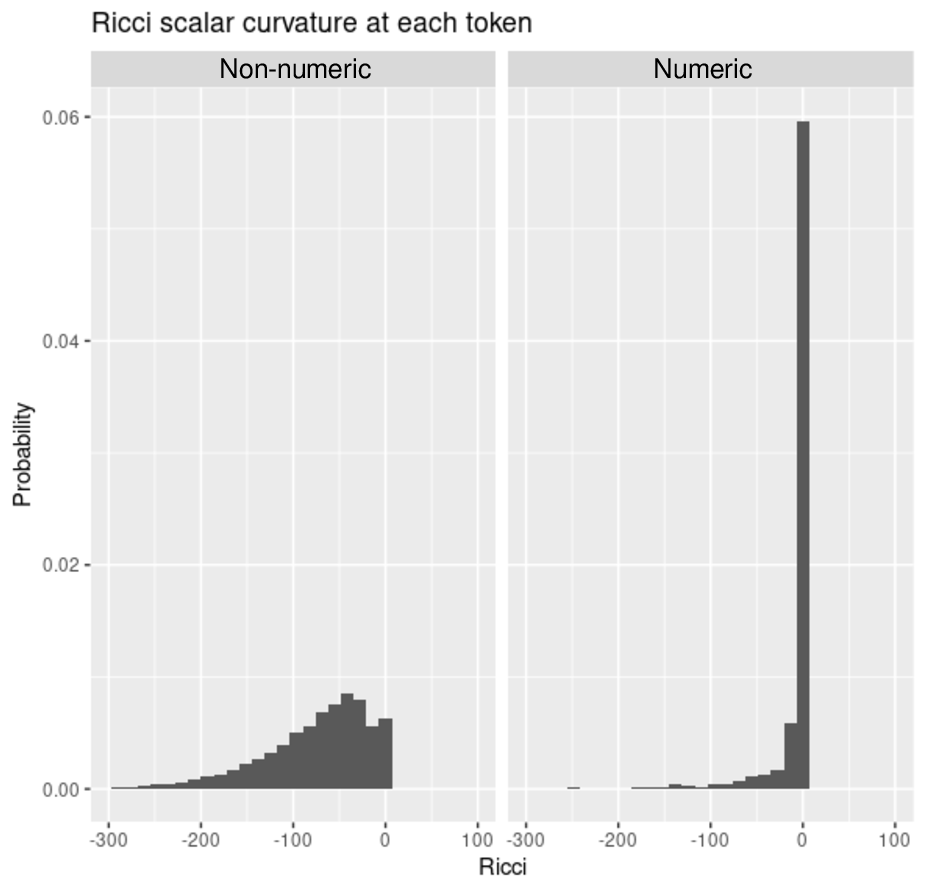}
    \caption{Histogram of estimated Ricci scalar curvature for GPT2.}
    \label{fig:gpt2_ricci_hist}
  \end{center}
\end{figure}

Figure \ref{fig:gpt2_ricci_hist} shows the distribution of estimates of Ricci scalar curvature for GPT2.  It is apparent that all estimates are negative, and the distribution of curvature estimates is unimodal.

\begin{figure}[!htbp]
  \begin{center}
    \includegraphics[width=4in]{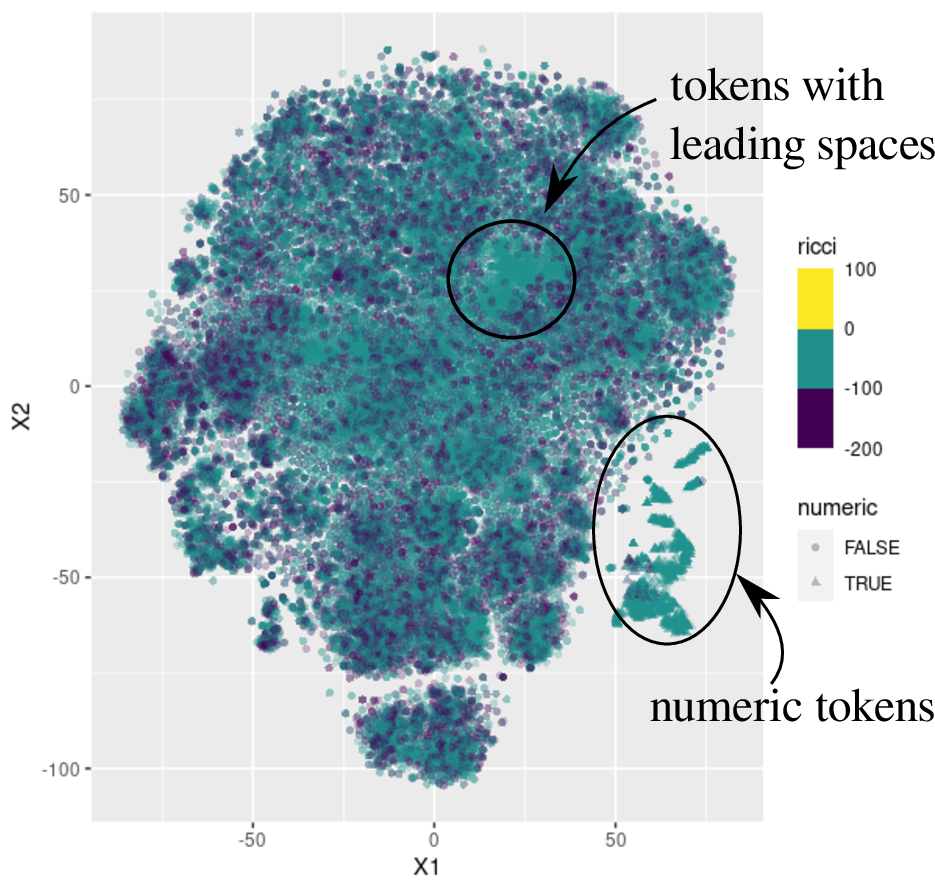}
    \caption{Estimated Ricci scalar curvature for GPT2, plotted using tSNE coordinates.   Curvature is indicated by the depth of color: darker points have higher (more negative) curvature.}
    \label{fig:gpt2_ricci_overview}
  \end{center}
\end{figure}

Figure \ref{fig:gpt2_ricci_overview} shows the Ricci curvature at each token, using the same tSNE coordinates as in Figure \ref{fig:gpt2_dim_overview}.
With one notable exception, the main connected component, consisting of mostly non-numeric tokens, has roughly constant negative Ricci scalar curvature.
The numeric tokens appear to lie on a region of substantially lower (closer to zero) curvature than the non-numeric tokens.
The uniformly low-curvature region within the non-numeric tokens appears to consist of tokens with leading spaces,
and may correspond to the separate low-dimensional stratum indicated along the slice in Figure \ref{fig:gpt2_dim_overview}.

\subsection{LLEMMA7B}

\begin{figure}[!htbp]
  \begin{center}
    \includegraphics[width=3in]{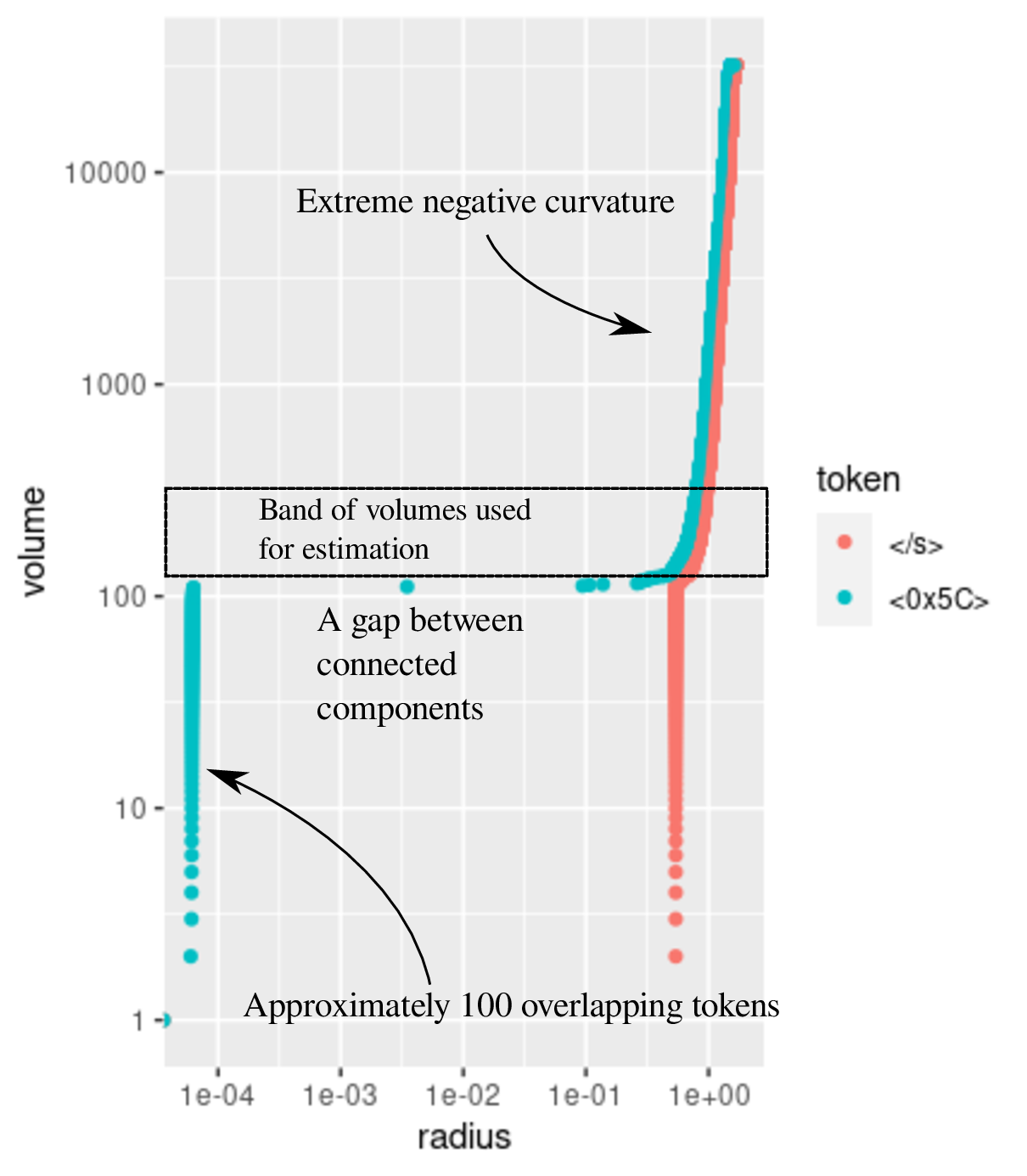}
    \caption{Volume-versus-radius for several tokens in LLEMMA7B.}
    \label{fig:llemma7b_v_r}
  \end{center}
\end{figure}

Figure \ref{fig:llemma7b_v_r} shows the volume-versus-radius curves for several tokens in LLEMMA7B.
The first thing to notice is that the nearest $100$ or so tokens are all equidistant from the two tokens shown.
This can be inferred by the vertical slopes on the left edge of both curves.
Since this is likely an artifact of how the embedding was designed, possibly due to rounding of the coordinates,
it is appropriate to exclude volumes below $100$ from our estimates below.
While clear stratifications are not visible as knees in the curves,
the token subspace has negative curvature near the two tokens shown since the curves are visibly concave up.

\begin{figure}[!htbp]
  \begin{center}
    \includegraphics[width=3in]{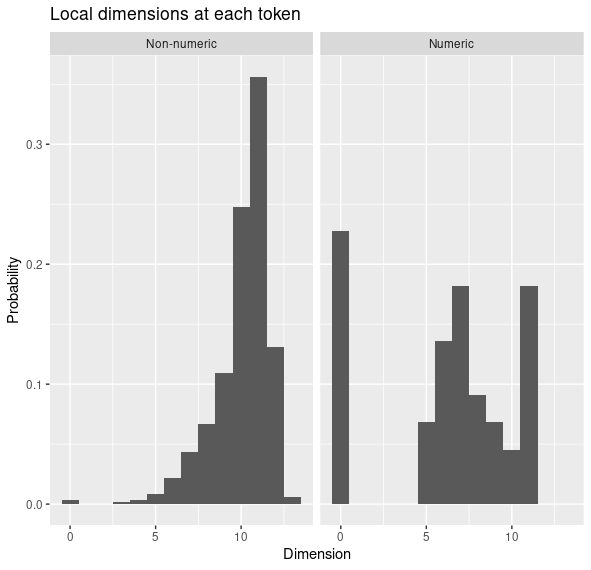}
    \caption{Histogram of estimated local dimensions for LLEMMA7B.  The large spike for numeric tokens of dimension $0$ corresponds to tokens that are isolated (not near any other tokens).}
    \label{fig:llemma7b_dim_hist}
  \end{center}
\end{figure}

Figure \ref{fig:llemma7b_dim_hist} shows the distribution of local dimensions for the token subspace of LLEMMA7B.
Again, the numeric tokens have a statistically significantly different distribution of local dimensions from the non-numeric tokens.
Indeed, many of the numeric tokens have a local dimension of $0$, which means that they are isolated points.
As such, many of the numeric tokens are not near each other.

\begin{figure}[!htbp]
  \begin{center}
    \includegraphics[width=4in]{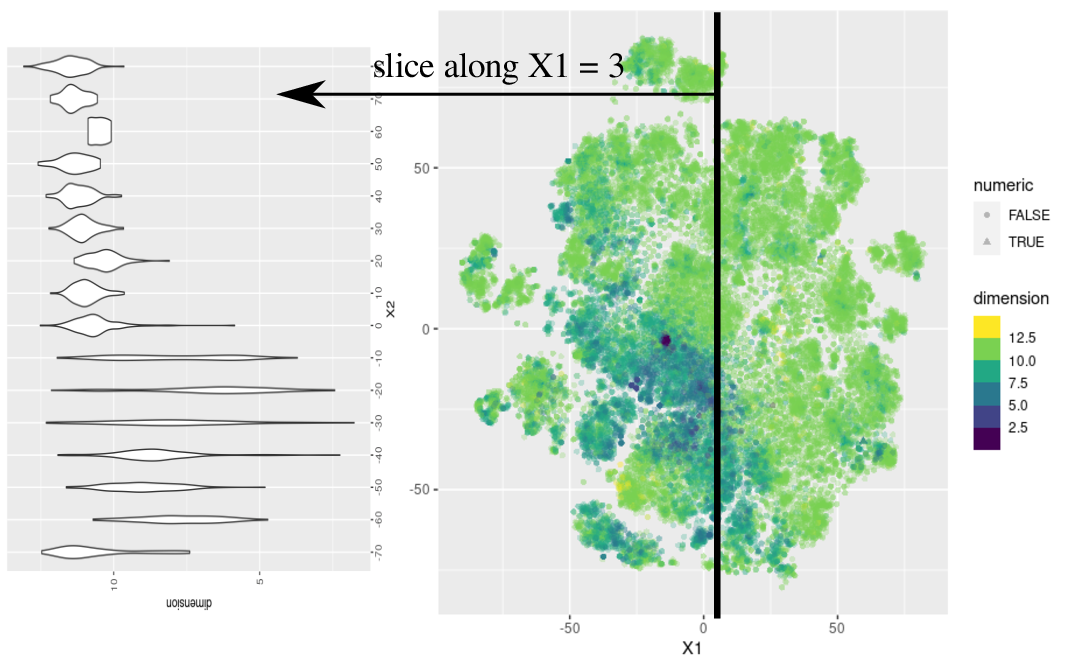}
    \caption{Estimated local dimension for LLEMMA7B, plotted using tSNE coordinates.  Dimension is indicated by the depth of color: darker points have lower local dimension.}
    \label{fig:llemma7b_dim_overview}
  \end{center}
\end{figure}

Unlike the case of GPT2, where the presence of knees in volume-versus-radius curves in Figure \ref{fig:llemma7b_v_r} establish that the token subspace is not a manifold, the stratifications in LLEMMA7B are a bit more subtle.
Figure \ref{fig:llemma7b_dim_overview} shows the resulting tSNE visualization, again colored by dimension.
The numeric tokens are located in two distinct regions of the plot as indicated.
There are several connected components visible in Figure \ref{fig:llemma7b_dim_overview}.
The largest connected component shows a marked transition in local dimension,
which is strong evidence of a stratification.

This stratification happens to be transverse to the $X1=3$ slice.
Figure \ref{fig:llemma7b_dim_overview} shows the local dimensions along that slice, which reveals that the local dimension is quite uniform on one side of the plot, and rather lower and more variable on the other.

Interestingly, there is a very small region in the center of the tSNE visualization in which the local dimension is quite small.
The tokens involved all appear to consist of non-printing Unicode characters.  Moreover, the same set of characters (and geometric features) is present in MISTRAL7B as well.

\begin{figure}[!htbp]
  \begin{center}
    \includegraphics[width=3in]{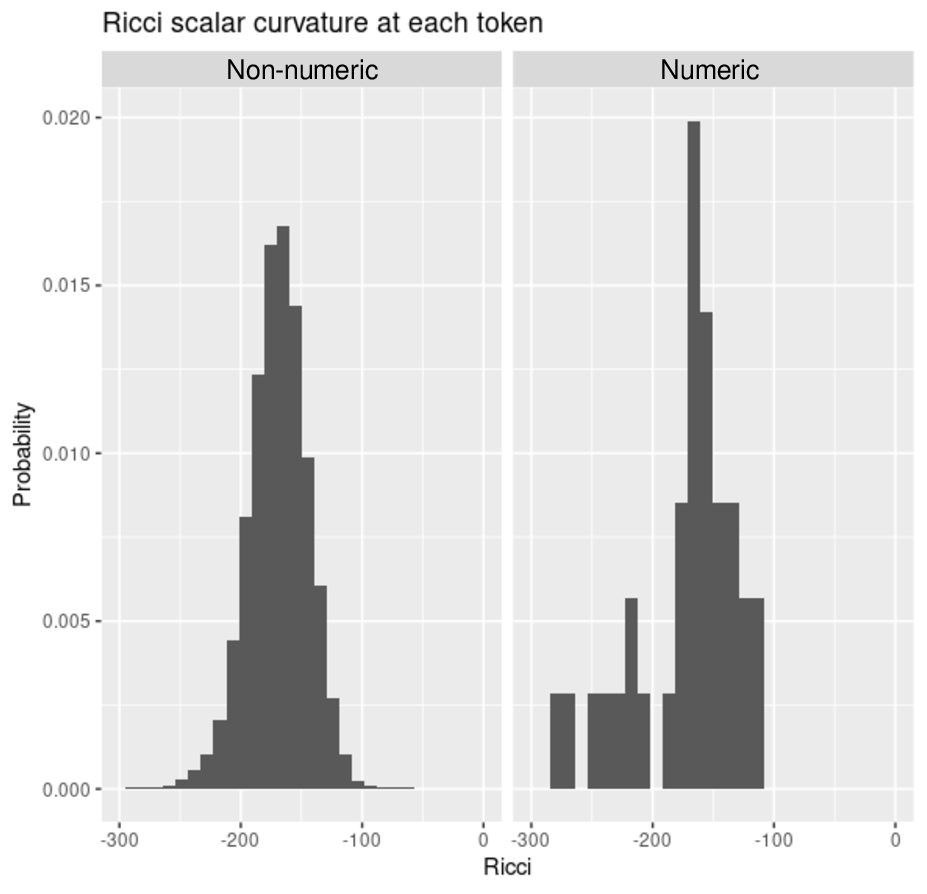}
    \caption{Histogram of estimated Ricci scalar curvature for LLEMMA7B.}
    \label{fig:llemma7b_ricci_hist}
  \end{center}
\end{figure}

\begin{figure}[!htbp]
  \begin{center}
    \includegraphics[width=4in]{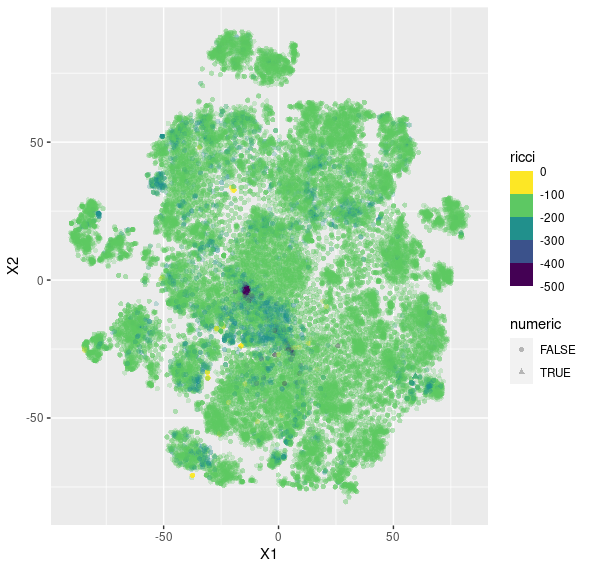}
    \caption{Estimated Ricci scalar curvature for LLEMMA7B, plotted using tSNE coordinates.  Curvature is indicated by the depth of color: darker points have higher (more negative) curvature.}
    \label{fig:llemma7b_ricci_overview}
  \end{center}
\end{figure}

Figure \ref{fig:llemma7b_ricci_hist} shows the histogram of Ricci scalar curvature for LLEMMA7B.
The distribution is  unimodal for both numeric and non-numeric tokens,
and is clearly negative.
Figure \ref{fig:llemma7b_ricci_overview} confirms this hypothesis: the Ricci scalar curvature is largely constant and negative,
though it changes abruptly along the stratification boundary.
It shows shows one clear exception to this pattern,
which is again the small low-dimensional region consisting of non-printing Unicode characters.

\subsection{MISTRAL7B}

\begin{figure}[!htbp]
  \begin{center}
    \includegraphics[width=3in]{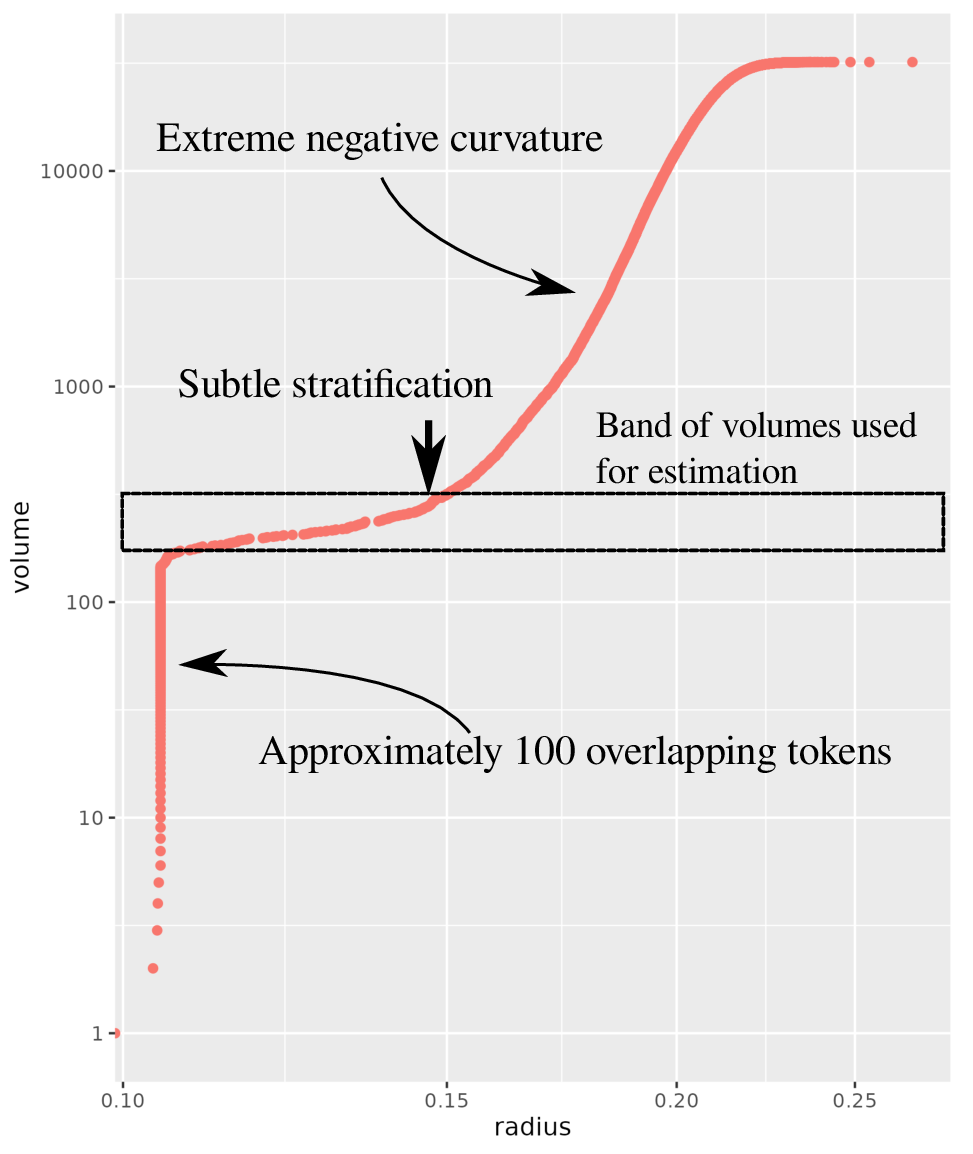}
    \caption{Volume-versus-radius for several tokens in MISTRAL7B.}
    \label{fig:mistral7b_v_r}
  \end{center}
\end{figure}

Figure \ref{fig:mistral7b_v_r} shows the volume-versus-radius for a typical token in MISTRAL.
As in the case of LLEMMA7B, there are about $150$ tokens that are equidistant from this token.
Since this appears to be the typical situation, we will use volumes with at least $150$ for estimation purposes below.
A knee in the curve corresponding to a stratification is visible, along with strong negative curvature.

\begin{figure}[!htbp]
  \begin{center}
    \includegraphics[width=3in]{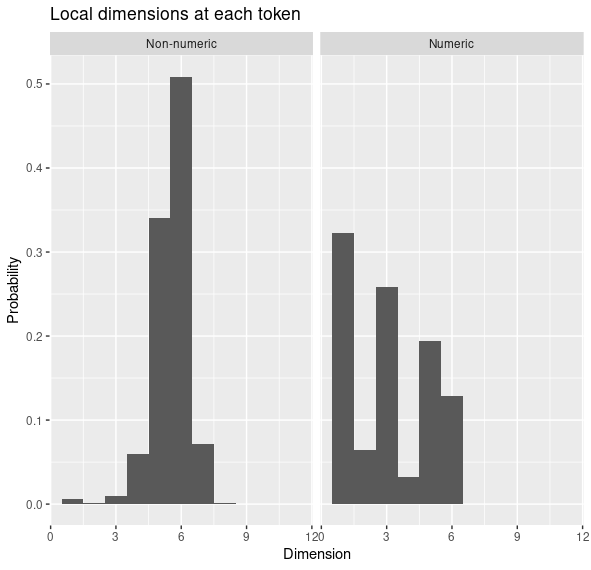}
    \caption{Histogram of estimated local dimensions for MISTRAL7B.}
    \label{fig:mistral7b_dim_hist}
  \end{center}
\end{figure}

Figure \ref{fig:mistral7b_dim_hist} shows the distribution of dimensions for all tokens in MISTRAL7B.
It is interesting to compare this with Figure \ref{fig:llemma7b_dim_hist}, because the same set of tokens is used in LLEMMA7B.
Both contain a large number of isolated (dimension $0$) numeric tokens.
Generally, the dimension of a given token is lower in MISTRAL7B than it is in LLEMMA7B.
This suggests that a given token in MISTRAL7B will tend to be easier for the model to distinguish than in LLEMMA7B,
because it has dramatically fewer neighbors.

\begin{figure}[!htbp]
  \begin{center}
    \includegraphics[width=4in]{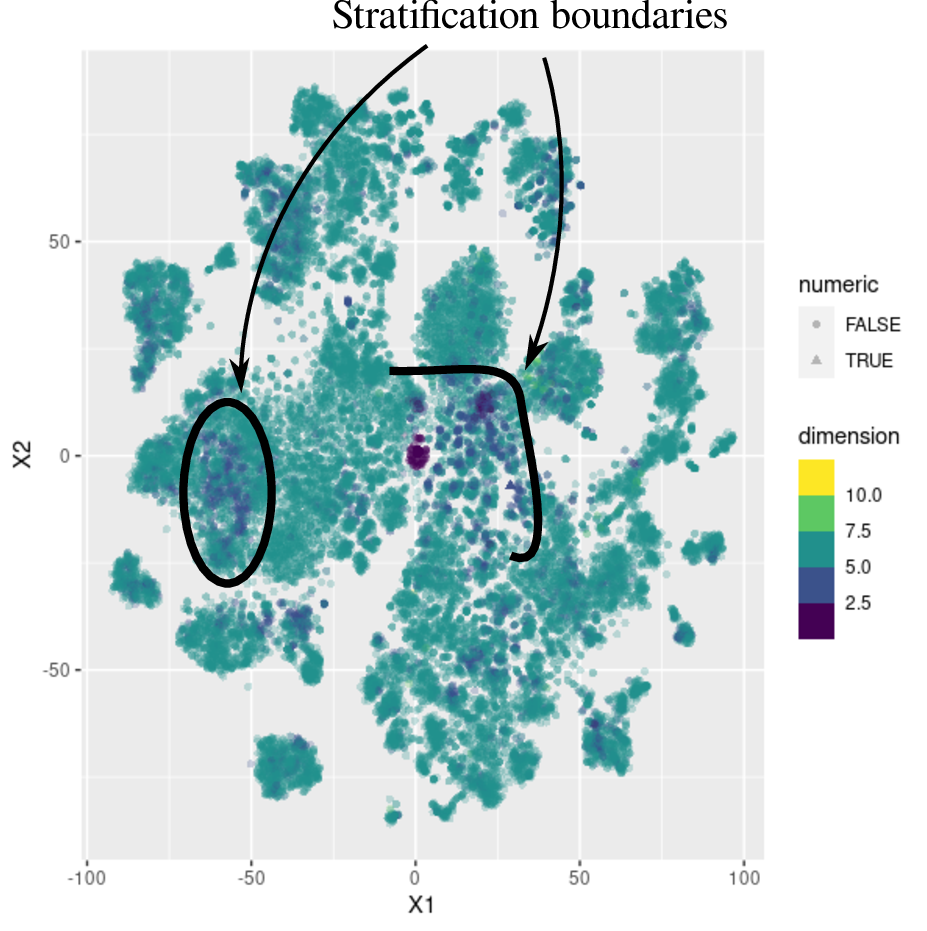}
    \caption{Estimated local dimension for MISTRAL7B, plotted using tSNE coordinates.  Dimension is indicated by the depth of color: darker points have lower local dimension.}
    \label{fig:mistral7b_dim_overview}
  \end{center}
\end{figure}

Figure \ref{fig:mistral7b_dim_overview} shows a tSNE embedding visualization of local dimension for MISTRAL7B.
Like the case of LLEMMA7B, there is a visible stratification boundary.
It also exhibits the same set of low-dimension non-printing Unicode characters,
which is probably unsurprising.
Unlike GPT2 and LLEMMA7B, the space for MISTRAL7B appears to be more disconnected.
However, due to the difficulties in interpreting clusters in tSNE embeddings, the reader is cautioned that this may not be representative.

\begin{figure}[!htbp]
  \begin{center}
    \includegraphics[width=3in]{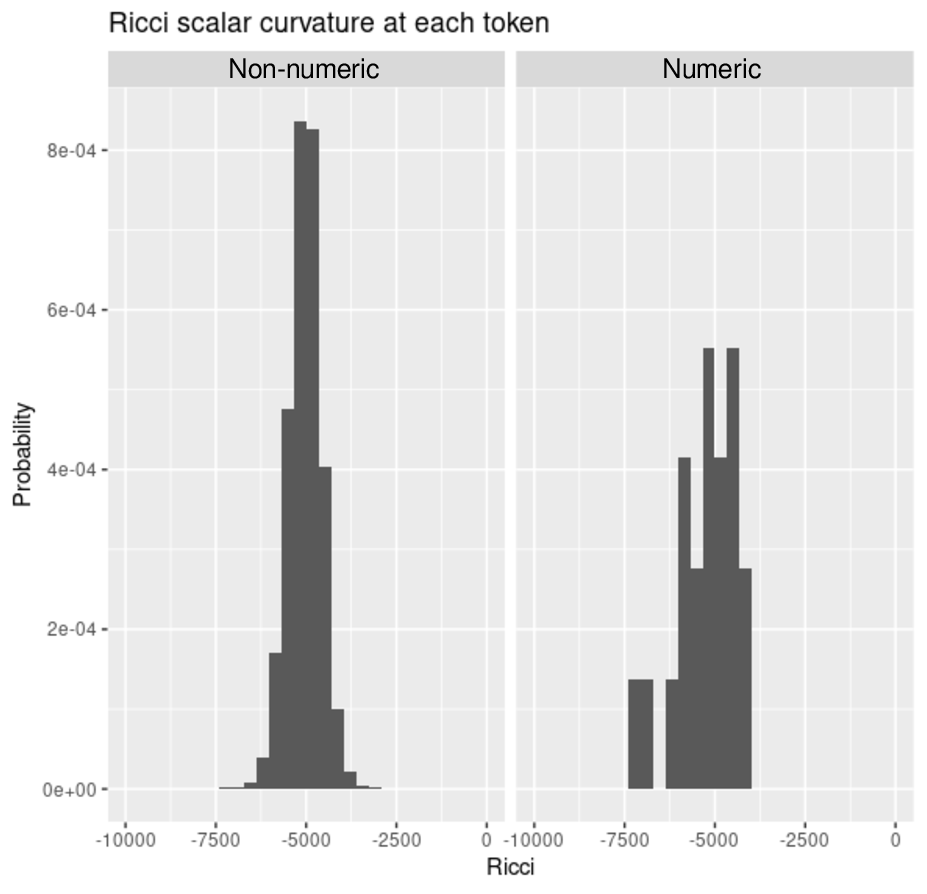}
    \caption{Histogram of estimated Ricci scalar curvature for MISTRAL7B.}
    \label{fig:mistral7b_ricci_hist}
  \end{center}
\end{figure}

\begin{figure}[!htbp]
  \begin{center}
    \includegraphics[width=4in]{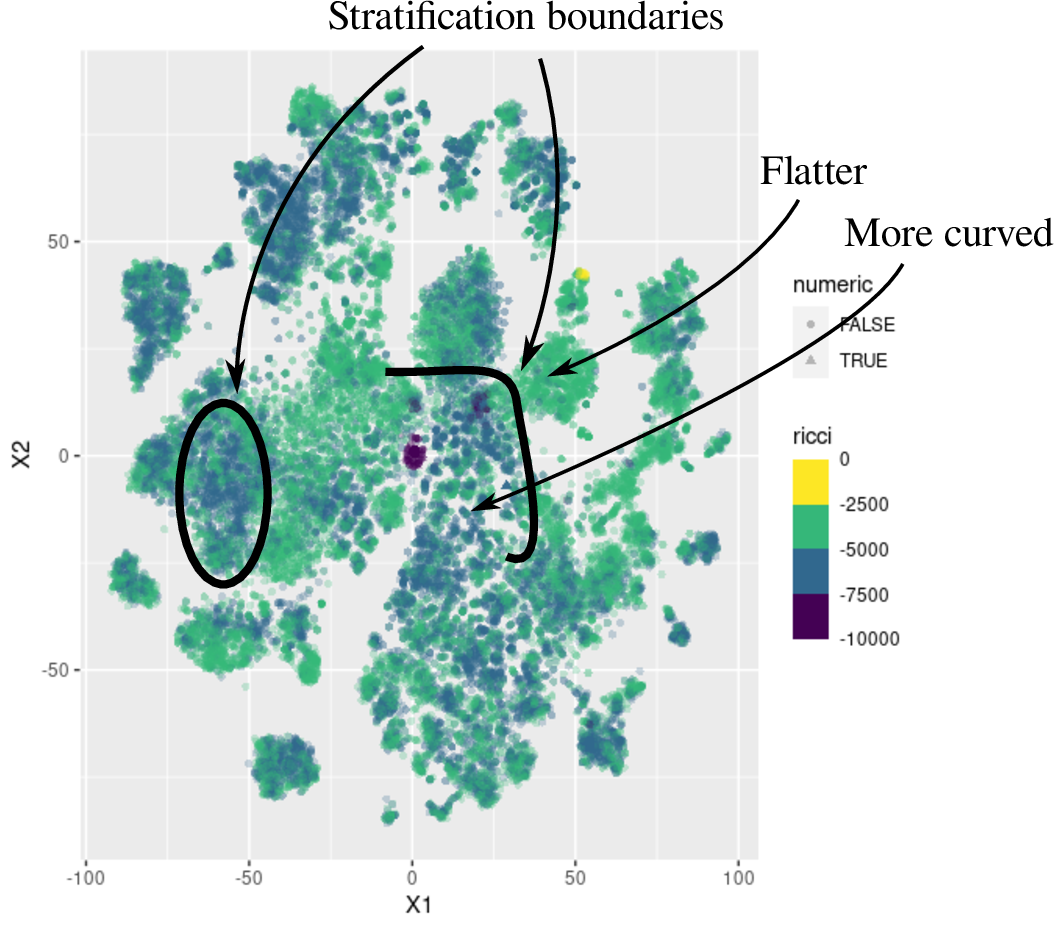}
    \caption{Estimated Ricci scalar curvature for MISTRAL7B, plotted using tSNE coordinates.  Curvature is indicated by the depth of color: darker points have higher (more negative) curvature.}
    \label{fig:mistral7b_ricci_overview}
  \end{center}
\end{figure}

Like the case of LLEMMA7B, the Ricci scalar curvature of the token subspace is mostly constant,
as indicated by the histograms in Figure \ref{fig:mistral7b_ricci_hist}.
This is confirmed in Figure \ref{fig:mistral7b_ricci_overview},
though the curvature is seen to vary somewhat over the token subspace.
Again, the collection of non-printing Unicode characters corresponds to a region of high curvature.

\section{Discussion}

We produced strong evidence that, for the three large language models under consideration in this article, the local dimension of the token subspace in large language models varies along connected components of that token subspace. This evidence implies that the token subspaces of these models cannot be manifolds.
Moreover, because the local dimension appears to change abruptly at multiple places along a slice through the token subspace, it seems likely that the token subspaces of these models consist of manifolds of different dimensions attached to each other. In other words, the token subspace of each of these three models are stratified manifolds.

Our work suggests that the presence and nature of the stratification of the token subspace impacts \emph{in a predictable way} the inferential quality and generative fluency of large language models. The examination of numeric and non-numeric tokens in this article, across models, provides an illustration of this point. While the token subspace of GPT2 is not a manifold, the \emph{numeric} token subspace is likely a manifold.
Since all of the numeric tokens are nearby each other and are confined to a constant dimension submanifold, this limits the expressivity of any continuous dynamic map (e.g. a transformer) that operates upon them. Because the dynamic map induced by a transformer is continuous (which follows from the continuity of the activation functions in a deep neural network \cite{villani_2024}),
this implies that any set of numeric tokens that are near each other will tend to be taken to other tokens that are also near each other.  It follows, simply from the finding that the numeric token subspace is (largely connected) manifold, that GPT2 will be a poor performer on mathematical queries: the numeric tokens of GPT2 are localized on a small number of connected components, such that that the model will tend to treat these numeric tokens as interchangeable, even though each token contains distinct and critically different semantic information. This insight is consistent with empirical observations of the \emph{inability} of GPT2 to consistently and coherently reason about numeric data. The topological and geometric structure of the token subspace of GPT2, a foundational model, is substantially different than the subspaces of both LLEMMA7B and MISTRAL7B, which were fine-tuned for performance on numeric and code-based queries. Both LLEMMA7B and MISTRAL7B have far fewer numeric tokens which are not located close to each other in the token subspace, consistent with these models' capability to distinguish between numeric tokens on inference.

Another implication of this work is how topological and geometric information together suggest model overfitting, which also has implications for generalizability and model performance against novel queries. In GPT2, the local dimension of the token subspace is approximately half the ambient latent space dimension; in both LLEMMA7B and MISTRAL7B, which are fine-tuned, the local dimension of the token subspace is much less than half the latent space dimension. One could infer the onset of overfitting behaviors based entirely upon the difference in dimension between the token subspace and the latent space,
which is usually called the \emph{codimension} of the subspace. Overfitting is a risk whenever the dimension of the latent space exceeds what is necessary to capture the structure of the token subspace.
According to the classical Whitney embedding theorem \cite{Lee_2003}, this kind of overfitting occurs whenever the codimension exceeds the token subspace dimension, and our findings show that this overfitting is definitely present in LLEMMA7B and MISTRAL7B, and is also present in portions of the token subspace for GPT2.

In addition to the issue of subspace codimension, we found that the token subspaces of all three models have strong negative Ricci curvature. Since the ambient latent space of each model is of zero curvature (because it is flat Euclidean space), this could explain why the latent space ``needs'' to have a much higher dimension than the token subspace.
Intuitively, a larger codimension for the token subspace (difference between the dimensions of the latent space and the token subspace) would give more ``room'' to accomodate the curved shape of the the token subspace.
However, according to classical interpolation error formulas \cite{isaacson1994analysis}, to which all smooth functions are subject,
the presence of high negative Ricci scalar curvature will dramatically increase the uncertainty of output of a function.
In the classical statistics literature, this is known as a consequence of overfitting.

An unfortunate anticipated consequence of this negative curvature is that retraining such as fine-tuning, which perturbs part of the token subspace locally, may result in precipitous global changes in the embedding. Large codimension in the presence of negative curvature means that linear transformations (which are usually present in transformers utilizing ReLU, among others \cite{villani2023any}),
will necessarily result in tokens being taken to points outside the token subspace. Local perturbations, as may be caused during fine-tuning, could result in bifurcations that yield global perturbations.
Some reflection of this has been observed in the literature \cite{Mishne_2022} because changing the latent space to a hyperbolic space (which has negative curvature) can improve trainging stability. Put another way, because of the stratifications and the negative curvature, bringing these points back to the token subspace is guaranteed to be unstable once they have strayed too far. In short, the geometry and topology of the subspace indicates that in portions of the token subspace, particularly with overfit, \emph{instability is mathematically unavoidable}, which translates to inferential uncertainty.

Taken together, our findings may explain why there appear to be instabilities in both the training and the use of large language models for query response, and moreover suggest that these instabilities are unavoidable in most cases.

\section*{Acknowledgements}

This article is based upon work supported by the Defense Advanced Research Projects Agency (DARPA).
Any opinions, findings and conclusions, or recommendations expressed in this material are those of the authors and do not necessarily reflect the views of DARPA.

\bibliographystyle{unsrt}
\bibliography{latentdimension_bib}

\section*{Appendix: Demonstration on spaces of known geometry and topology}

To illustrate the viability of this method beyond its application to the token subspaces of large language models,
we demonstrate that it works on several known manifolds with different metrics.
Although the method does not presuppose a manifold, it is easy to compute the correct parameters directly on circles, spheres, and disks.
The results of these experiments are shown in Table \ref{tab:space_fittings}.

\begin{table}[!htpb]
  \caption{Estimated parameters of known manifolds}
  \label{tab:space_fittings}
  \begin{tabular}{|l|c|c|c||c|c|c|}
    \hline
    Manifold & Metric & Parameter & True & \multicolumn{3}{|c|}{Estimated} \\
    &&  & value  & Q1 & Q2 & Q3 \\
    \hline
    \hline
    Circle & Arclength & Dimension & $1$ & $0.945$ & $$0.978$$ & $1.01$ \\
    & $R=1$ & Scaling & $2$ & $1.91$ & $1.96$ & $2.16$ \\    
    &  & Ricci & $0$ & $-0.175$ & $-0.0181$ & $0.0849$ \\
    \hline
    Circle & Euclidean & Dimension & $1$ & $0.910$ & $0.964$ & $1.03$ \\
    & $R=1$ & Scaling & $2$ & $1.84$ & $1.96$ & $2.06$ \\    
    &  & Ricci & $-3/4$ & $-2.05$ & $-1.29$ & $-0.492$ \\
    \hline

    Sphere & Arclength & Dimension & $2$ & $1.89$ & $1.97$ & $2.06$ \\
    & $R=1$ & Scaling & $\pi$ & $3.11$ & $3.18$ & $3.24$ \\    
    &  & Ricci & $2$ & $1.71$ & $2.04$ & $2.41$ \\
    \hline

    Sphere & Euclidean & Dimension & $2$ & $1.75$ & $1.90$ & $2.09$ \\
    &  $R=1$ & Scaling & $\pi$ & $2.74$ & $3.02$ & $3.37$ \\    
    &  & Ricci & $0$ & $-2.87$ & $-0.906$ & $1.05$ \\
    \hline

    Disk & Euclidean & Dimension & $2$ & $1.78$ & $1.91$ & $2.03$ \\
    & $R=1/2$ & Scaling & $\pi$ & $1.57$ & $2.45$ & $3.42$ \\
    &  & Ricci & $0$ & $-65$ & $15.7$ & $117$ \\
    \hline
  \end{tabular}
\end{table}

We report quartiles (rather than mean and standard deviation) for the estimated values because all of their respective distributions are highly non-normal.
(The $p$-values from the Shapiro-Wilk test for normality were far below $0.01$ in every case tested, which supports this claim.)
The dimension, scaling coefficient, and Ricci curvature are defined by Equation \eqref{eq:volume_curvature} above.
We were able to determine the scaling coefficient directly because the volume of each of these spaces is known in advance.
Since we tested many points on each of these spaces, we show the quartiles for each parameter computed by our method.
The interquartile range of the distribution of estimates of each example contains the correct answer, which indicates that our method is producing the correct answers on average.
Derivations of the parameters from classical formulas are included below.
In addition to the manifolds listed in Table \ref{tab:space_fittings},
we also ran the method against a stratified manifold consisting of three strata: a $2$-dimensional disk joined to a $1$-dimensional circle at a $0$-dimensional point. 

\begin{figure}[!htbp]
  \begin{center}
    \includegraphics[width=3in]{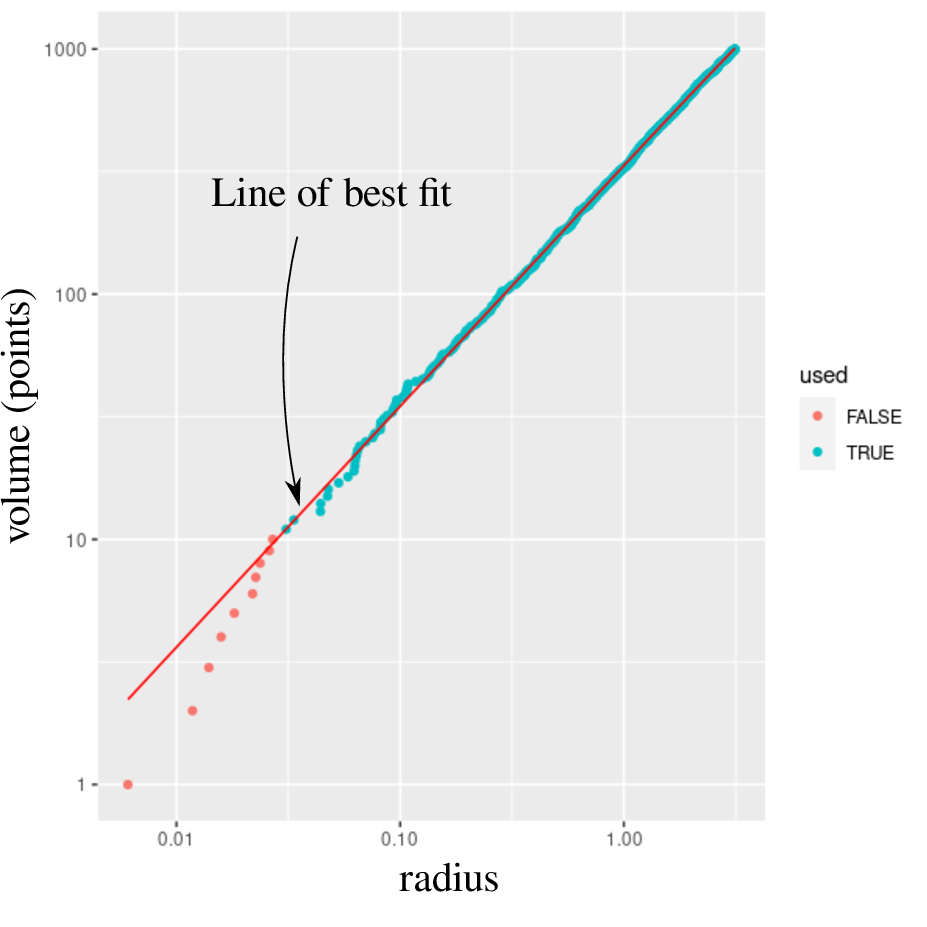}
    \caption{Log-log regression problem for a random point on the circle of radius $1$.  Points used in the regression are marked in green.  The line of best fit from the data as well as the true volume dependence are both shown.}
    \label{fig:circle_ll_estimate}
  \end{center}
\end{figure}

As a first example, consider the circle of radius $1$, in which the distance between points is given by the arclength between them.
The total ``volume'' of the circle is therefore its circumference, namely $2 \pi$.
Figure \ref{fig:circle_ll_estimate} shows the regression problem specified in Equation \eqref{eq:volume_log_matrix} for a typical point on the circle.
Notice that the slope of the line of best fit is close to the ideal one, which has slope exactly $1$.
Table \ref{tab:space_fittings} indicates that the dimension is within the interquartile range of all of the dimension estimates.
In addition to using the arclength distance, we can measure the distance between points on the circle embedded in the plane by using the Euclidean distance on the plane.
Because this is not a Riemannian metric, the Ricci scalar curvature is spuriously nonzero, which our method correctly detects.

\begin{figure}[!htbp]
  \begin{center}
    \includegraphics[width=3in]{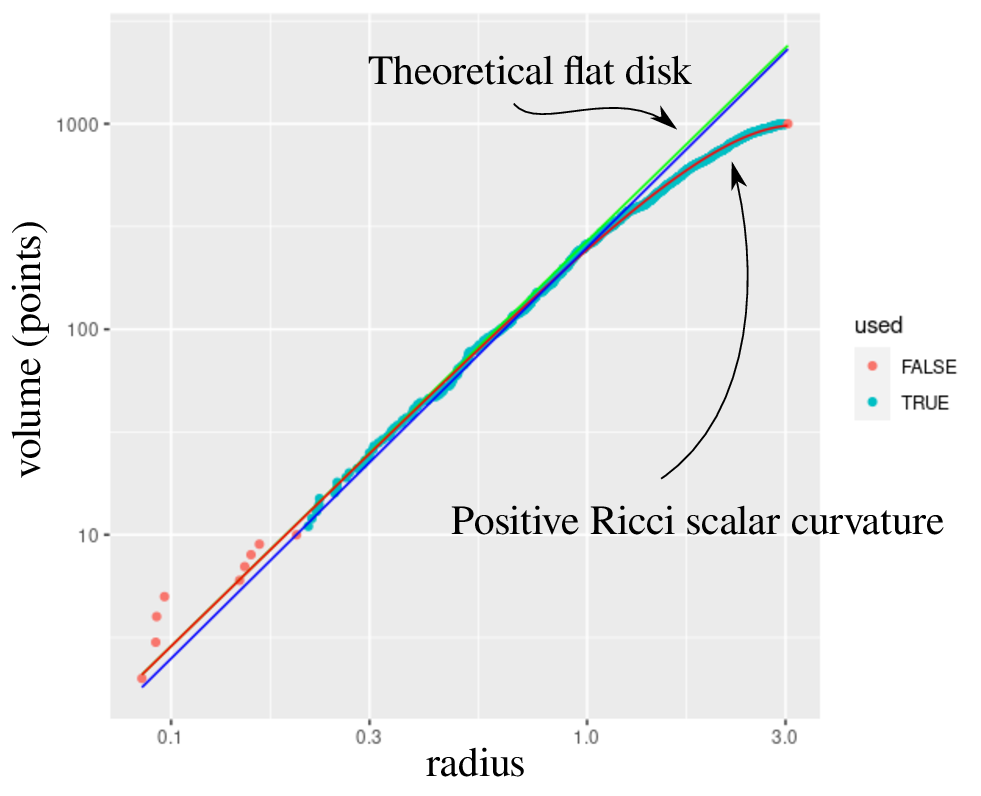} 
    \caption{Log-log regression problem for a random point on the $2$-dimensional sphere of radius $1$ with the arclength metric.  Points used in the regression are marked in green.  Notice that the curve is concave down, which indicates that the Ricci scalar curvature is positive according to Equation \ref{eq:volume_curvature_log}.}
    \label{fig:sphere_ll_estimate}
  \end{center}
\end{figure}

The same two metrics can be used for a $2$-dimensional sphere embedded isometrically in $\mathbb{R}^3$.
Again, Table \ref{tab:space_fittings} shows that the correct values of the parameters are estimated by our method.
Notice that points on the sphere are determined by two values, latitude and longitude, so the dimension of the sphere is $2$.

In the case of the sphere with the arclength metric, the curvature is positive (see Table \ref{tab:space_fittings} and the Appendix).
According to Equation \ref{eq:volume_curvature_log}, this implies that the volume of a disk on the sphere should be somewhat less for a given radius than for a disk on a flat space.
This is confirmed in Figure \ref{fig:sphere_ll_estimate} because the curve is concave down, because according to the discussion in Section \ref{sec:v_r} this corresponds to positive Ricci scalar curvature.

\begin{figure}[!htbp]
  \begin{center}
    \includegraphics[width=4in]{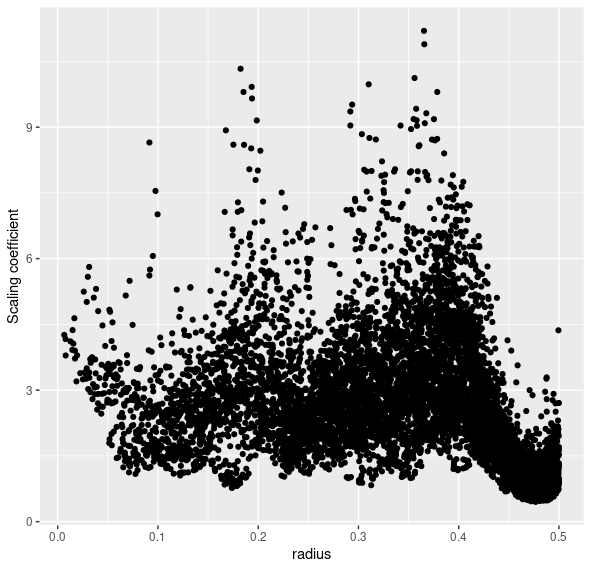}
    \caption{Effect of the boundary on the scaling coefficient for a disk of radius $R=0.5$.}
    \label{fig:disk_boundary_effect}
  \end{center}
\end{figure}

For the case of a flat disk of radius $1/2$ isometrically embedded in the plane,
the dimension is $2$.
The scaling coefficient will be different for points within the interior versus points on the boundary.
At points away from the boundary on the interior,
\begin{equation*}
  v = \pi r^2,
\end{equation*}
so scaling coefficient is $\pi$
but for points on the boundary,
\begin{equation*}
  v < \frac{\pi}{2} r^2,
\end{equation*}
so the scaling coefficient is less than $\pi/2$.
Figure \ref{fig:disk_boundary_effect} shows the estimates of scaling coefficient as a function of radius.
For small radii, the scaling coefficient is indeed around $\pi$, but this drops sharply as the boundary is approached.
Because our method uses distances to neighboring points to pose the regression problem,
the boundary does impact points near, but not on, the boundary.

\begin{figure}[!htbp]
  \begin{center}
    \includegraphics[width=3in]{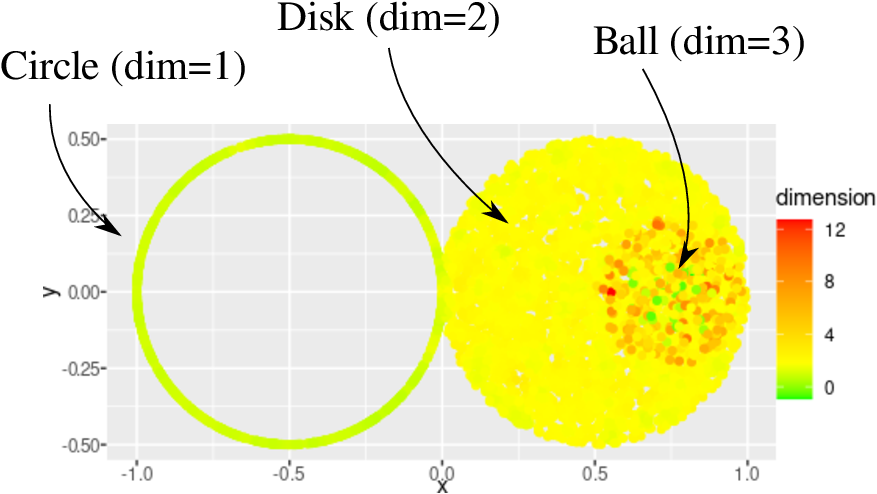}
    \caption{Estimated local dimension for a stratified space.  The circular arc has dimension $1$, the disk has dimension $2$, and the ball has dimension $3$.}
    \label{fig:stratified_overview}
  \end{center}
\end{figure}

As a final example, consider a stratified manifold consisting of the union of a disk and a circle joined at a single point, and a ball attached to part of the disk as shown in Figure \ref{fig:stratified_overview}.

\begin{figure}[!htbp]
  \begin{center}
    \includegraphics[width=3in]{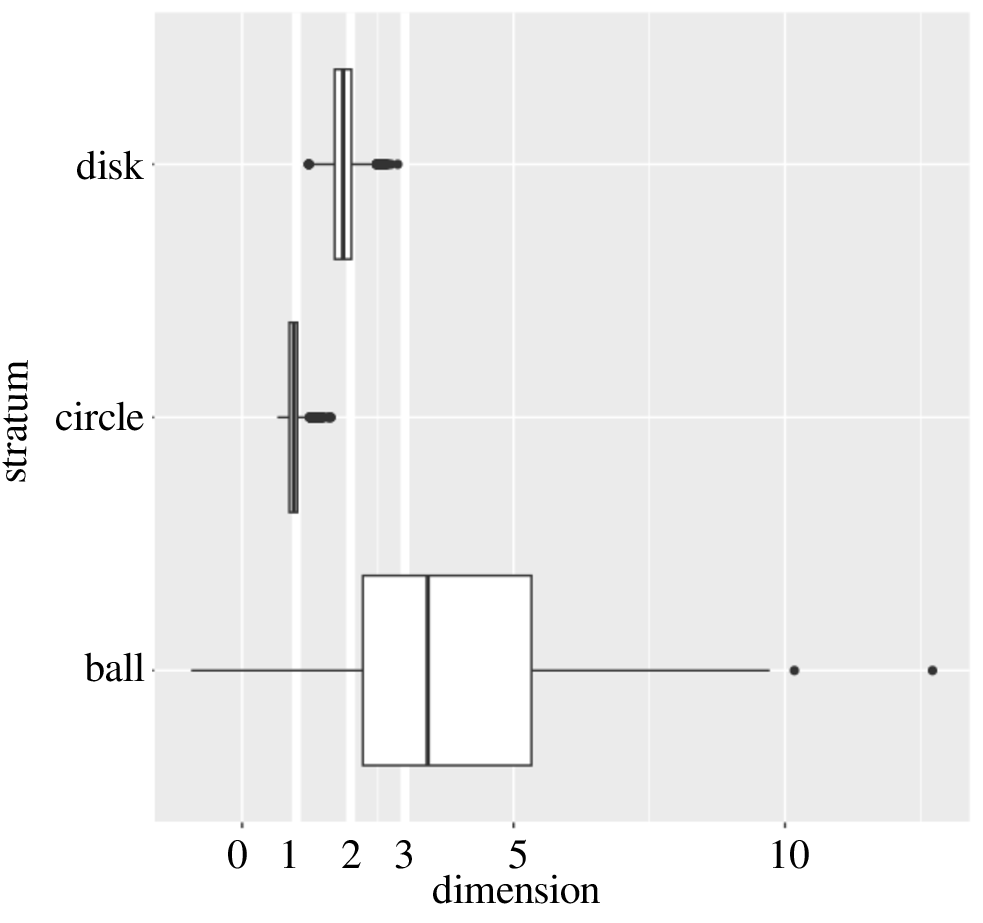}
    \caption{Distribution of local dimension estimates for a stratified space.  The circular arc has dimension $1$, the disk has dimension $2$, and the ball has dimension $3$.} 
    \label{fig:stratified_dimension}
  \end{center}
\end{figure}

The points within the interior of the ball have local dimension $3$, those on the interior of the disk have local dimension $2$, and the points on the circle have local dimension $1$.
This is correctly estimated by our method, as shown in Figure \ref{fig:stratified_dimension}.

\begin{figure}[!htbp]
  \begin{center}
    \includegraphics[width=3in]{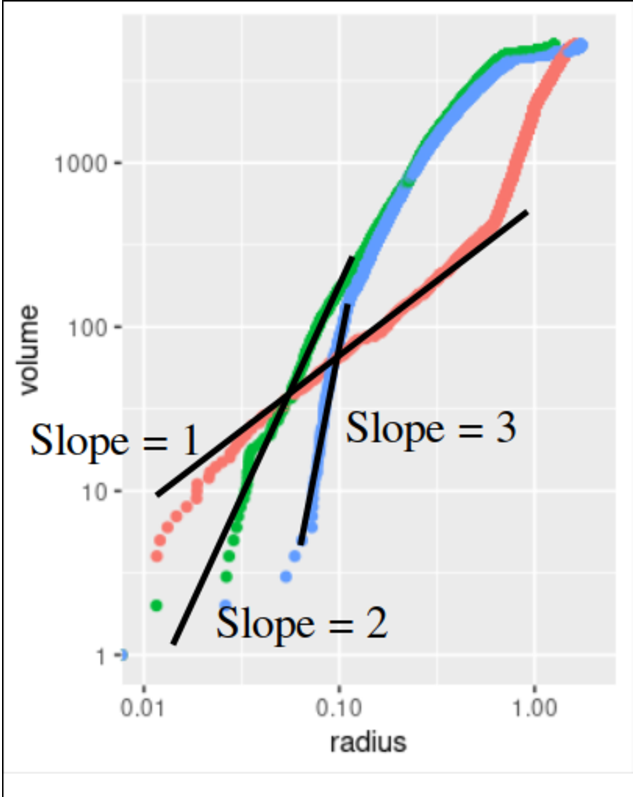}
    \caption{Volume-versus-radius for the stratified space shown in Figure \ref{fig:stratified_overview}.}
    \label{fig:stratified_v_r}
  \end{center}
\end{figure}

As explained in Section \ref{sec:v_r}, volume-versus-radius curves with knees imply the presence of detect stratifications.
These knees are prominently visible in Figure \ref{fig:stratified_v_r}, which confirms the presence of stratifications.

Of the manifolds listed in Table \ref{tab:space_fittings}, the circle and sphere with the Euclidean metric are \emph{not} Riemannian manifolds, so Equation \eqref{eq:volume_curvature} does not apply.
This is easily seen in the case of the circle: the distance between antipodes is twice the radius, but this disagrees with the arclength distance, namely $\pi$ times the radius.
Nevertheless, we can still use Equation \eqref{eq:volume_curvature} to estimate the Ricci scalar curvature from the data \emph{as if the metric were Riemannian}.
The justification for this is that the leading term is still correct and the remaining terms are geometrically meaningful, though not necessarily interpreted as Ricci scalar curvature \cite[Chap. 5]{Yomdin_2004}.
Doing this, we will obtain a topological manifold with the same dimension and scaling coefficient, but the curvature will be changed as is clear from the entries in Table \ref{tab:space_fittings}.

\begin{figure}[!htbp]
  \begin{center}
    \includegraphics[width=2.5in]{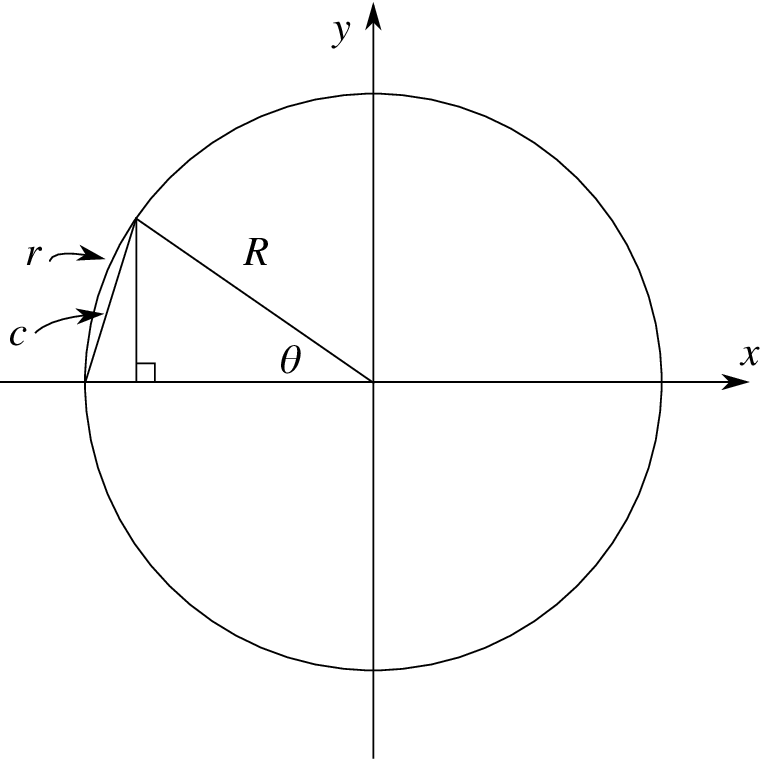}
    \caption{Geometry for a sector of both a circle in the plane and a planar cross section of sphere, showing its radius $R$, arclength $r$, chord length $c$, and included angle $\theta$.}
    \label{fig:circle_setup}
  \end{center}
\end{figure}

Consider Figure \ref{fig:circle_setup}, showing a sector of a circle (or cross section of a sphere) of radius $R$.
Since both the circle and the sphere are homogeneous spaces, we may choose any point as a center.
It is most convenient to select the point furthest on the left of the diagram.
We wish to compute the volume of a ball centered on this point in four different ways:
\begin{enumerate}
\item For the circle with its usual metric corresponding to arclength.  Because the circle has one free parameter, it is $1$-dimensional, so ``volume'' means arclength.
\item For the circle with distances between points being measured using the Euclidean metric in the plane.  Because the circle has one free parameter, it is $1$-dimensional, so ``volume'' means arclength.
\item For the sphere with its usual metric corresponding to great circle distance.  Because the sphere has two free parameters (latitude and longitude), it is $2$-dimensional, so ``volume'' means area.
\item For the sphere with distances between points being measured using the Euclidean metric in the $\mathbb{R}^3$.  Because the sphere has two free parameters (latitude and longitude), it is $2$-dimensional, so ``volume'' means area.
\end{enumerate}

The two cases of metrics on the circle can be summarized by the following,
\begin{equation}
  \label{eq:circle}
  \begin{aligned}
  v &= \begin{cases}
    2r & \text{using arclength}\\
    2r\left(1+\frac{r^2}{24R^2} + O(r^4)\right) & \text{using Euclidean distance}
  \end{cases}\\
  &=  2r + O(r^2).
  \end{aligned}
\end{equation}
For case (1), the circle with arclength, since radius and arclength are exactly the same, the only thing to recognize is that the ball centered at a given point is symmetric about that point.
Hence the ``volume'' of the ball is twice its radius, in accordance with Equation \eqref{eq:circle}.

The two cases of metrics on the sphere are summarized by
\begin{equation}
  \label{eq:sphere}
  \begin{aligned}
  v &= \begin{cases}
    \pi r^2\left(1 - \frac{r^2}{12R^2} + O(r^4)\right) & \text{using arclength}\\
    \pi r^2  & \text{using Euclidean distance}
  \end{cases}\\
  &=  \pi r^2 + O(r^2).
  \end{aligned}
\end{equation}
It is a little surprising that the case of Euclidean distance on the sphere yields the familiar formula for the area of a disk on the plane,
but this is due to the fact that the Euclidean distance \emph{is not} a Riemannian metric for the sphere.
Nevertheless, both of these expressions are derived below.

In both the circle and the sphere, the length of the chord $c$ depends on the angle $\theta$, which in turn depends on arclength via $r = R\theta$.
Therefore,
\begin{equation}
  \label{eq:chordal}
  \begin{aligned}
    c &= \sqrt{(R-R\cos \theta)^2 + R^2 \sin \theta}\\
    &= \sqrt{2} R \sqrt{1 - \cos \theta} \\
    &= \sqrt{2} R \sqrt{1 - \cos (r/R)} \\
  \end{aligned}
\end{equation}
Since the ``volume'' of the ball is twice the arclength, we need to invert the above expression, namely
\begin{equation*}
  \begin{aligned}
    v &= 2 r = 2 R \arccos \left(1 - \frac{c^2}{2R^2}\right) \\
    &= 2R \left( \frac{c}{R} + \frac{c^3}{24 R^3} + O(c^5) \right) \\
     &= 2c  \left(1 + \frac{c^2}{24 R^2} + O(c^4)\right),\\
  \end{aligned}
\end{equation*}
which agrees with the second case of Equation \eqref{eq:circle} after noting that the chord $c$ is the Euclidean distance between points.

For the case of the sphere, the ``volume'' is the area of the surface of revolution swept by the arc with length $r$.
Recognizing that points on the circle are defined by
\begin{equation*}
  x^2 + y^2 = R^2,
\end{equation*}
we can use the classical formula for surface area of revolution,
\begin{equation*}
  \begin{aligned}
    A(r) &= 2\pi \int_{-R}^{-R\cos (r/R)} y(x) \sqrt{1 + y'(x)^2} dx \\
    &= 2\pi \int_{-R}^{-R\cos (r/R)} \sqrt{R^2 - x^2} \sqrt{1 + \frac{x^2}{R^2 - x^2}} dx\\
    &= 2\pi \int_{-R}^{-R\cos (r/R)} R \; dx = 2\pi R^2 \left(1-\cos \frac{r}{R}\right).
  \end{aligned}
\end{equation*}
Combining the above expression with Equation \eqref{eq:chordal} results in the expression for surface area as a function of the Euclidean distance,
\begin{equation*}
  A(c) = \pi \left(\sqrt{2} R \sqrt{1 - \cos (r/R)}\right)^2 = \pi c^2,
\end{equation*}
which is the second case of Equation \eqref{eq:sphere}.
On the other hand, a Taylor expansion for the volume in the case of the arclength distance yields the first case, namely
\begin{equation*}
  \begin{aligned}
    A(r) &= 2\pi R^2 \left(1-\cos \frac{r}{R}\right) \\
    &= 2 \pi R^2 \left(1 - \left(1 - \frac{1}{2R^2}r^2 + \frac{1}{24R^4}r^4 + O(r^6)\right)\right) \\
    &= \pi r^2 \left(1 - \frac{1}{12 R^2}r^2 + O(r^4)\right).
  \end{aligned}
\end{equation*}

\end{document}